\documentclass{amsart}
\usepackage{amsfonts}
\input{epsf}
\newtheorem{theorem}{Theorem}[section]

\theoremstyle{definition}

\theoremstyle{remark}

\numberwithin{equation}{section}

\begin{document}

\title{Properties of the wealth process in a market microstructure model}

\author{Ted Theodosopoulos}
\address{Department of Decision Sciences and \\ Department of Mathematics 
\\ Drexel University \\ Philadelphia, PA 19066}
\email{theo@drexel.edu}
\urladdr{www.lebow.drexel.edu/theodosopoulos}
\author{Ming Yuen}
\address{Department of Decision Sciences 
\\ Drexel University \\ Philadelphia, PA 19066}
\email{mmy23@drexel.edu}


\date{February 1, 2005.}


\begin{abstract}
In this short paper we define the wealth process in a spin model for market microstructure, for individual agents and in aggregate.  The agents in our model try to balance their desire to belong to the local majority (herding behavior), defined over random network neighborhoods, and the occasional advantage of belonging to the global minority (contrarian trading).  We arrive at a classification of the martingale properties of this wealth process and use it to determine the strategic stability of the agents' interactions.  Our goal is to add a behavioral interpretation to this stochastic agent-based model for market fluctuations.  
\end{abstract}

\maketitle

\section{Introduction}

Over the past two years a series of models have been proposed that attempt to study the microstructure of financial transactions.  An increasing portion of these models are motivated from analogy to physical systems that exhibit the type of distributed decision making and frustrated coordination regimes that characterize market empirically \cite{bouchaud4}.  We are interested in adding some economic motivation to the statistical mechanics of a class of spin market models.

The class of models we will investigate aims to capture the trade-off between two competing drives that motivate market participants.  On the one hand, each agent in the interaction network receives input (e.g. advice, information, opinions) from a randomly chosen neighborhood.  The receiving agents tries to conform their behavior to the norm expressed in this local interaction.  On the other hand, each agent receives information about the global imbalance between buy and sell interest in the market.  Under certain conditions, it is preferable for the agent to position themselves in the global minority, in anticipation of a move by the entire network.  

The former consideration is the basis of the traditional voter process \cite{granovsky, majka}, which is used to model herding behavior.  The latter has given rise to the minority game and its variants \cite{burgos, challet}, as models for contrarian trading.  Bonrholdt and collaborators combined these two incentives into the Hamiltonian for a spin system that models market microstructure \cite{bornholdt1}.  A variant of this process was studied analytically in \cite{theo1} and the invariant measure was characterized.  Also, an uncertainty relation was observed, which poses a fundamental limit to the amount of control one can have on a market of this type.

Here was use the representation of the process arrived at in \cite{theo1} to define and investigate the wealth process for individual agents and the market as a whole.  We start by reviewing the model and establishing our notation.  We then proceed to present our main result, which leads to a classification of the possible qualitative behaviors for the wealth process.  In the following sections we characterize the strategic stability of the wealth process, and analyze the residual risk for the wealth of an individual agent and the aggregate wealth of the market.  Finally we show some numerical results regarding the path properties of the stochastic process that go beyond the invariant measure from \cite{theo1}.  

\section{Model Description}

The model we use is a variant of Bornholdt's spin market model \cite{bornholdt1} that was presented in \cite{theo1}.  The state space $X$ of our model is the set of spin configurations on a lattice on the $d$-dimensional torus\footnote{Here we use the notation ${\mathcal T}^d$ to denote the object $\underbrace{ {\mathcal S}^1 \times \ldots \times {\mathcal S}^1}_d$.} $Y \doteq \left({\mathcal Z}/L \right)^d \subset {\mathcal T}^d$, i.e. $X \subset \{-1,1 \}^Y$, for an appropriately chosen $L$ so that $|Y|=N$.  The path of a typical element of $X$ is given by $\eta: Y \times ( 0, \infty ) \longrightarrow \{-1,1\}$ and each site $x \in Y$ is endowed with a (typically $\ell_1$) neighborhood ${\mathcal N} (x) \subset Y$ it inherits from the natural topology on the torus ${\mathcal T}^d$.  In this paper a version of 'rapid stirring' \cite{durrett} is applied by randomizing the neighborhood structure ${\mathcal N}(\cdot)$.  In particular, for each $x \in Y$, ${\mathcal N} (x)$ is a uniformly chosen random subset of $Y$, of cardinality $2d$.  To be more precise, for a set $A$ and a positive integer $k$, let $F(A,k) = \left\{ \left(a_1, a_2, \ldots, a_k \right) \in A^k | a_i \neq a_j \mbox{         } \forall i,j=1,2,\ldots,k \right\}$.  Then, for any $x \in Y$, let $\left\{ {\mathcal N} (x, \cdot) \right\}$ be a family of iid uniform random variables taking values in $F \left( Y \setminus \{x\}, 2d \right)$.

We construct a continuous time Markov process with transitions occurring at exponentially distributed epochs, $T_n$, with rate 1.  We proceed to construct a transition matrix for the spins, based on the following interaction potential:
$$h(x,T_n) = \sum_{y \in {\mathcal N}(x,n)} \eta(y,T_n) - \alpha \eta(x,T_n) N^{-1} \left|\sum_{y \in Y} \eta(y,T_n) \right|,$$
where $\alpha>0$ is the coupling constant between local and global interactions.  At time $T_n$ (i.e. the $n$th epoch) a random site $x$ is chosen and its spin is changed to $+1$ with probability $p^+ \doteq \left( 1+ \exp \left\{- 2\beta h \left(x,T_n \right) \right\} \right)^{-1}$ and to $-1$ with probability $p^- = 1- p^+$, where $\beta$ is the normalized inverse temperature.  We define the price and volume processes as follows:
\begin{eqnarray*}
P(t) & = & P^\ast (t) \exp \left\{ \lambda N^{-1} \sum_{y \in Y} \eta(y,t) \right\} \\ 
V(t) & = & N^+ (t) \vee N^- (t)
\end{eqnarray*}
where $a \vee b = \max\{a,b\}$, $a \wedge b = \min\{a,b\}$, $\lambda$ is a parameter, $P^\ast$ is an exogenous previsible `fundamental' price process \cite{bornholdt1} (which we assume identically equal to 1 for purposes of the current study) and $N^\pm (t) \doteq \left| \left\{y \in Y \/ | \/ \eta(y,t)=\pm1 \right\} \right|$ denotes the number of sites with a positive or negative spin respectively.
Using the auxiliary variables $\bar{X}_n = N^+ \left( T_n \right) - N^+ \left( T_{n-1} \right)$, $X_n = \left|\bar{X}_n \right|$ and $Y_n = \left| 2N^+ \left(T_n \right) - N \right|$, we can express $V \left(T_n \right) = \left( Y_n + N \right)/2$ and the volatility, 
\begin{equation}
\sigma \left( \left. \log {\frac {P \left(T_n \right)}{P \left( T_{n-1} \right)}} \right| {\mathcal F}_{n-1} \right) = {\frac {2 \lambda}{N}} \sqrt{{\rm P} \left( X_n = 0 \right) {\rm P} \left( X_n = 1 \right)}. \label{eq:pricevoly}
\end{equation}
The invariant measure of the state variable $N^+$ was characterized in \cite{theo1}.  Specifically, it was shown that when the coupling constant $\alpha$ is large compared to the network connectivity (specifically greater than the degree, $2d$) then the process is an the \textit{supercritical regime}, which maintains persistent market excitation even in the \textit{frozen} phase of this spin model ($\beta \rightarrow \infty$).  

\section{Martingale analysis of the wealth process}

In the current paper we restrict our attention to the frozen phase of the supercritical regime of this market microstructure model.  We proceed to define the wealth process of agent $y$ as 
$$W(y,T_k) \doteq \left\{
\begin{array}{ll}
W(y,T_{k-1}) & \mbox{for $J_{k}=y$;} \\
W(y,T_{k-1}) + \eta (y,T_{k-1})\Delta P(T_k) & \mbox{otherwise.}\\
\end{array}
\right.,$$
where $J_k \in \{1, 2, \ldots, N\}$ is the agent chosen on the $k^{\rm th}$ epoch $T_k$, and $W(y,0)= K(y)$ representing the initial capital available to agent $y$.  Our interest will be in characterizing the martingale properties of $W(y, \cdot)$ as well as the aggregate wealth, $W$, defined as  
\begin{eqnarray}
W(T_k)&=& W(T_{k-1})+\Delta P(T_k) \sum_{y \in Y} \eta (y,T_{k-1}) \nonumber \\
& = & W(T_{k-1})+\Delta P(T_k) \left[2N^+(T_{k-1})-N+\bar{X}_k \right]. \label{eq:globalwealth}
\end{eqnarray}
It will be helpful in our presentation to develop the following notation, where $i$ represents the state $N^+$ of the Markov process:
\begin{eqnarray*}
& & c = \left\lceil \alpha \left| {\frac {i}{N}} - {\frac {1}{2}} \right| \right\rceil , \;\;\;\;\; f_+(i,j) = C_{\scriptscriptstyle j}^{\scriptscriptstyle i-1} C_{\scriptscriptstyle 2d-j}^{\scriptscriptstyle N-i} \left/ C_{\scriptscriptstyle 2d}^{\scriptscriptstyle N-1} \right.  \;\;\;\;\; f_-(i,j) = C_{\scriptscriptstyle j}^{\scriptscriptstyle i} C_{\scriptscriptstyle 2d-j}^{\scriptscriptstyle N-i-1} \left/ C_{\scriptscriptstyle 2d}^{\scriptscriptstyle N-1} \right. \\
& & \begin{array}{cc} P_{\scriptscriptstyle ++}(i) = P_{\scriptscriptstyle --}(i) = 0 & \;\;\;\;\; \mbox{if $i \not\in \left[N \left( {\frac {1}{2}} - {\frac {d}{\alpha}} \right), N \left( {\frac {1}{2}} + {\frac {d}{\alpha}} \right) \right]$} \end{array} \\
& & \begin{array}{cc} P_{\scriptscriptstyle ++}(i) = {\frac {i}{N}} {\displaystyle \sum_{\scriptscriptstyle j=(d + c) \vee (i+2d-N)}^{\scriptscriptstyle 2d \wedge i}} f_+(i,j) & \;\;\;\;\; \mbox{if $i \in \left[N \left( {\frac {1}{2}} - {\frac {d}{\alpha}} \right), N \left( {\frac {1}{2}} + {\frac {d}{\alpha}} \right) \right]$} \end{array} \\
& & \begin{array}{cc} P_{\scriptscriptstyle --}(i) = \left(1 -{\frac {i}{N}} \right) {\displaystyle \sum_{\scriptscriptstyle j=0 \vee (i+2d-N)}^{\scriptscriptstyle (d-c) \wedge i}} f_-(i,j) & \;\;\;\;\; \mbox{ if $i \in \left[N \left( {\frac {1}{2}} - {\frac {d}{\alpha}} \right), N \left( {\frac {1}{2}} + {\frac {d}{\alpha}} \right) \right]$} \end{array} \\
& & P_{\scriptscriptstyle +-} = {\frac {i}{N}} -P_{\scriptscriptstyle ++} , \;\;\;\;\; P_{\scriptscriptstyle -+} = \left(1 -{\frac {i}{N}} \right) -P_{\scriptscriptstyle --}
\end{eqnarray*}
where $C_k^n$ denotes the combinations $n$ choose $k$.

We begin by observing that
\begin{equation}
\Delta P(T_k) = P(T_{k-1}) \left\{e^{\frac {2 \lambda \bar{X}_k}{N}}-1 \right\}. \label{eq:deltaP}
\end{equation}
Therefore
\begin{eqnarray*}
{\bf E} \left[ e^{\frac {2 \lambda \bar{X}_k}{N}}-1 \left| {\mathcal B}_{k-1} \right. \right] & = & \left( e^{\frac {2 \lambda}{N}}-1 \right) P_{\scriptscriptstyle +-} + \left( e^{\frac {2 \lambda}{N}}-1 \right) P_{\scriptscriptstyle -+} \\
& = & \left( e^{\frac {2 \lambda}{N}}-1 \right) \left(P_{\scriptscriptstyle -+} -e^{-{\frac {2 \lambda}{N}}} P_{\scriptscriptstyle +-} \right),
\end{eqnarray*}
and
\begin{eqnarray*}
{\bf E} \left[\bar{X}_k \left( e^{\frac {2 \lambda \bar{X}_k}{N}}-1 \right) \left| {\mathcal B}_{k-1} \right. \right] & = & -\left( e^{\frac {2 \lambda}{N}}-1 \right) P_{\scriptscriptstyle +-} + \left( e^{\frac {2 \lambda}{N}}-1 \right) P_{\scriptscriptstyle -+} \\
& = & \left( e^{\frac {2 \lambda}{N}}-1 \right) \left(P_{\scriptscriptstyle -+} +e^{-{\frac {2 \lambda}{N}}} P_{\scriptscriptstyle +-} \right),
\end{eqnarray*}
Putting these observations together and using (\ref{eq:globalwealth}) and (\ref{eq:deltaP}) we see that
\begin{eqnarray*}
{\bf E} \left[W (T_k) -W(T_{k-1}) \left| {\mathcal B}_{k-1} \right. \right] & = & P(T_{k-1}) \left\{ \left[ 2N^+ (T_{k-1}) -N \right] {\bf E} \left[e^{\frac {2 \lambda \bar{X}_k}{N}}-1 \left| {\mathcal B}_{k-1} \right. \right] + \right. \\
& & \left. + {\bf E} \left[\bar{X}_k \left( e^{\frac {2 \lambda \bar{X}_k}{N}}-1 \right) \left| {\mathcal B}_{k-1} \right. \right] \right\} =
\end{eqnarray*}
\begin{equation}
= P(T_{k-1}) \left( e^{\frac {2 \lambda}{N}}-1 \right) \left\{ \left[2N^+ (T_{k-1}) -N+1 \right] P_{\scriptscriptstyle -+} -e^{-{\frac {2 \lambda}{N}}}\left[2N^+ (T_{k-1}) -N-1 \right] P_{\scriptscriptstyle +-} \right\} \label{eq:globalmart}
\end{equation}

\begin{figure}
\epsfxsize=6in
\epsfbox{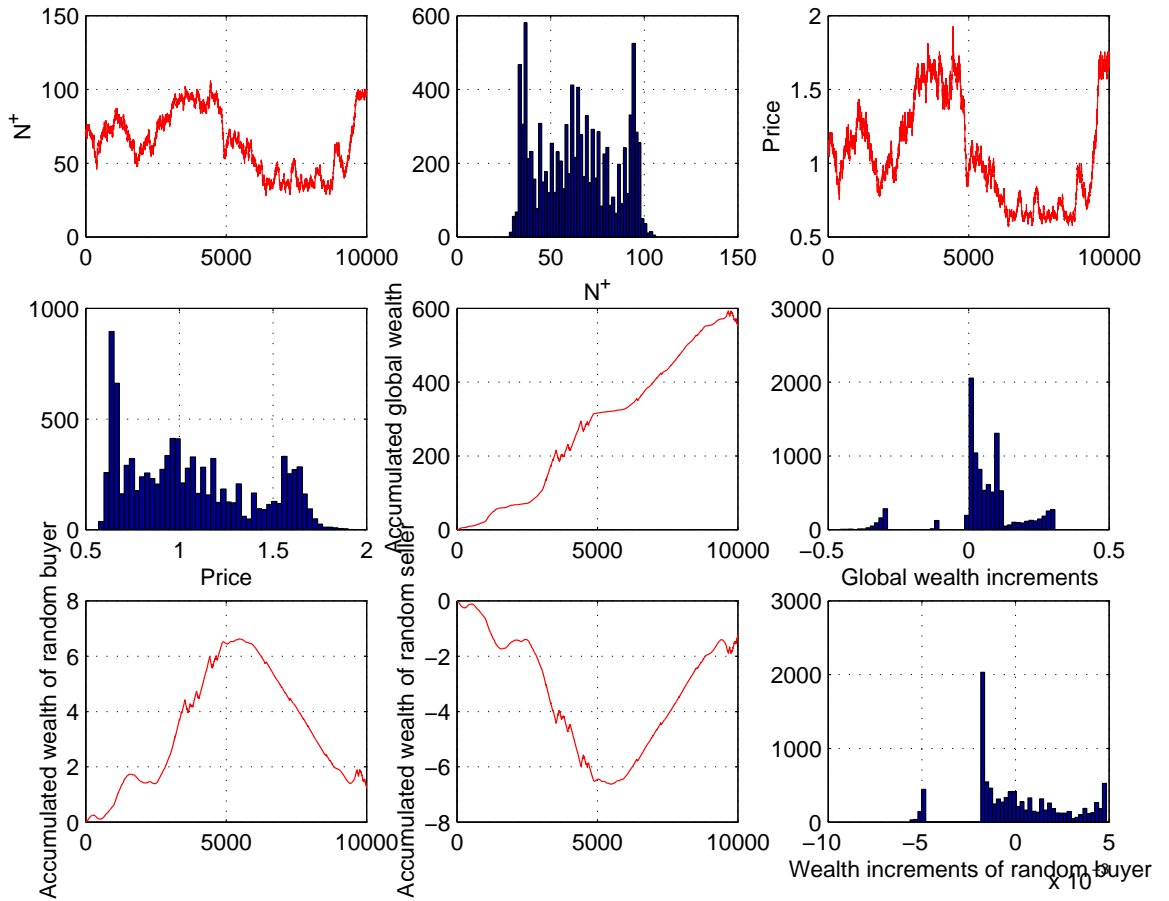}
\caption{A simulation of the wealth process in the frozen phase of the supercritical regime for 10,000 steps ($d=2, N=128, \alpha=4.1, \lambda=1$).  The last three frames are included for illustration purposes only.  They assume that the agent maintains their spin throughout the simulation, which is clearly a very unlikely event.}
\label{fig:pathsimulation}
\end{figure}

A similar analysis leads to the following statement for the conditional expectation of the increments of the wealth for agent $y$:
\begin{eqnarray}
& & {\bf E} \left[W (y, T_k) -W(y, T_{k-1}) \left| {\mathcal B}_{k-1} \right. \right] = \nonumber \\
& & = \eta(y, T_{k-1}) P (T_{k-1}) \left(1 -{\frac {1}{N}} \right) \left( e^{\frac {2 \lambda}{N}}-1 \right) \left\{ P_{\scriptscriptstyle -+} -e^{-{\frac {2 \lambda}{N}}} P_{\scriptscriptstyle +-} \right\} \label{eq:localmart}
\end{eqnarray}
Figure \ref{fig:pathsimulation} shows a sample path of $N^+$, $W$.  As expected, the distribution of $N^+$ converges to the trimodal invariant measure described in \cite{theo1}.  
The following theorem uses the conditional expectations of the wealth increments arrived at in (\ref{eq:globalmart}) to characterize the stochastic dynamics of the wealth process:
\begin{theorem}
\label{theorem1}
There exist four integer functions $g_1(\alpha)$, $g_2(\lambda)$, $g_3(\lambda)$ and $g_4(\alpha, \lambda)$ such that the aggregate wealth process $W$ is a submartingale \cite{feller} while the Markov process $N^+$ is the region $\left[g_1(\alpha), g_2(\lambda) \right] \cup \left[ g_3(\lambda), g_4 (\alpha,\lambda) \right]$.  Moreover, $g_1$ is a increasing function of $\alpha$, $g_2$ is a decreasing while $g_3$ is an increasing function of $\lambda$ and $g_4$ is increasing with respect to $\lambda$ and decreasing with respect to $\alpha$.  Two qualitative transitions occur as $\alpha$ and $\lambda$ vary:
\begin{enumerate}
\item For large enough $\lambda$, $g_2$ becomes equal to $g_1$, and then the lower interval disappears.
\item On the other hand, for small enough $\lambda$, $g_2(\lambda)$ becomes equal to $g_3(\lambda)$ and the two intervals merge, while $g_1(\alpha)= \arg \sup_{i<N/2} \pi_\infty (i)$ 
and $\lim_{\lambda \rightarrow 0} g_4(\alpha, \lambda)= \arg \sup_{i>N/2} \pi_\infty (i)$, where $\pi_\infty$ is the invariant measure of the Markov process $N^+$ in the frozen phase \cite{theo1}.
\end{enumerate}
\end{theorem}
The table below illustrates the functions $g_1, g_2, g_3, g_4$ and their dependence on $\alpha$ and $\lambda$.

\hspace{1in}
\vspace{0.2in}
\begin{tabular}{|c|c|c|c|c|c|c|c|} \hline
$N$ & $\alpha$ & $d$ & $\lambda$ & $g_1$ & $g_2$ & $g_3$ & $g_4$ \\ \hline
128 & 5 & 2 & 7 & 39 & & & 89 \\
128 & 5 & 2 & 8 & 39 & 44 & 54 & 89 \\
128 & 5 & 2 & 9 & 39 & 40 & 56 & 89 \\
128 & 5 & 2 & 10 & & & 57 & 89 \\
128 & 5 & 2 & 64 & & & 63 & 109 \\
128 & 4.1 & 2 & 7 & 33 & & & 95 \\
128 & 4.1 & 2 & 8 & 33 & 44 & 54 & 95 \\
128 & 4.1 & 2 & 9 & 33 & 40 & 56 & 95 \\
128 & 4.1 & 2 & 10 & 33 & 38 & 57 & 95 \\
128 & 4.1 & 2 & 11 & 33 & 36 & 58 & 95 \\
128 & 4.1 & 2 & 12 & 33 & 34 & 59 & 95 \\
128 & 4.1 & 2 & 13 & 33 & 33 & 59 & 95 \\
128 & 4.1 & 2 & 14 & & & 60 & 95 \\
128 & 4.1 & 2 & 64 & & & 63 & 109 \\
\hline
\end{tabular}
\vspace{0.2in}

\section{Strategic stability and path properties}

We saw earlier in (\ref{eq:pricevoly}) that the parameter $\lambda$ is proportional to the price volatility in this spin market model.  Thus one can conclude from Theorem \ref{theorem1} that sufficiently volatile markets experience a submartingale aggregate wealth process (i.e. the market generates wealth on average) if and only if there is a surplus of buyers.

Using (\ref{eq:localmart}) and the same logic as above we can deduce that the individual wealth process follows a qualitatively similar process.  The primary difference is that the individual wealth process increments depend on the spin at the site under consideration.  Thus, when the aggregate wealth follows a submartingale with a surplus of sellers (which, as we saw before, can only occur for low enough price volatility) the buyers follow a supermartingale.  If they were allowed to flip their spins independently, rather than follow the stochastic process, they would choose to do so, in order to experience a submartingale wealth process.  In so doing they would reduce $N^+$, bringing closer to $g_1$.  

\begin{figure}
\epsfxsize=4in
\epsfbox{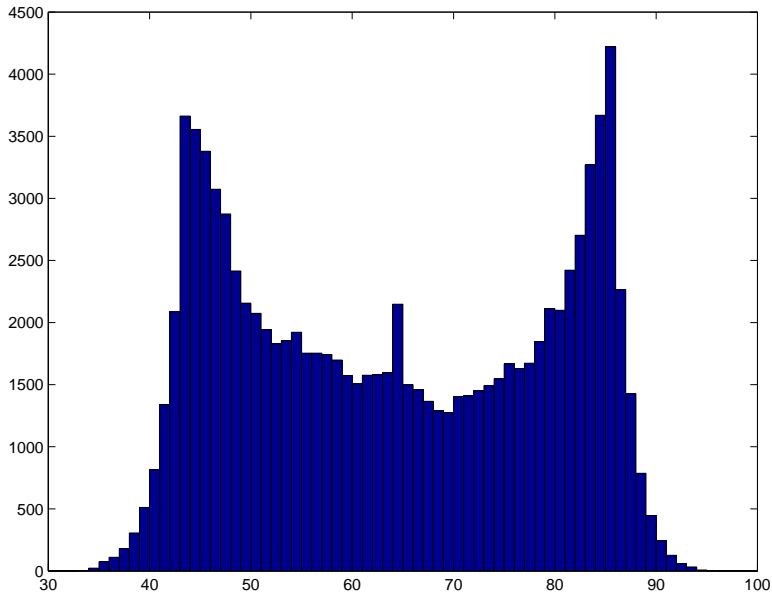}
\caption{This histogram shows the amount of time the Markov process spends in the state $N^+ = x$ as a function of $x$.  The simulation from which this histogram is constructed ran for $10^5$ steps with $d=2, N=128, \alpha=6, \lambda=1$.}
\label{fig:invariant}
\end{figure}

The same strategic argument can be applied to a market with a surplus of buyers.  As a result, of this argument and Theorem \ref{theorem1}, we can deduce that these stochastic market dynamics have two strategically stable equilibria, at $g_1$ and $g_4$.  We would therefore expect the market to spend most of the time in the vicinity of these two states.  Figure \ref{fig:invariant} confirms this expectation for a random sample path for the process.

\begin{figure}
\epsfxsize=6in
\epsfbox{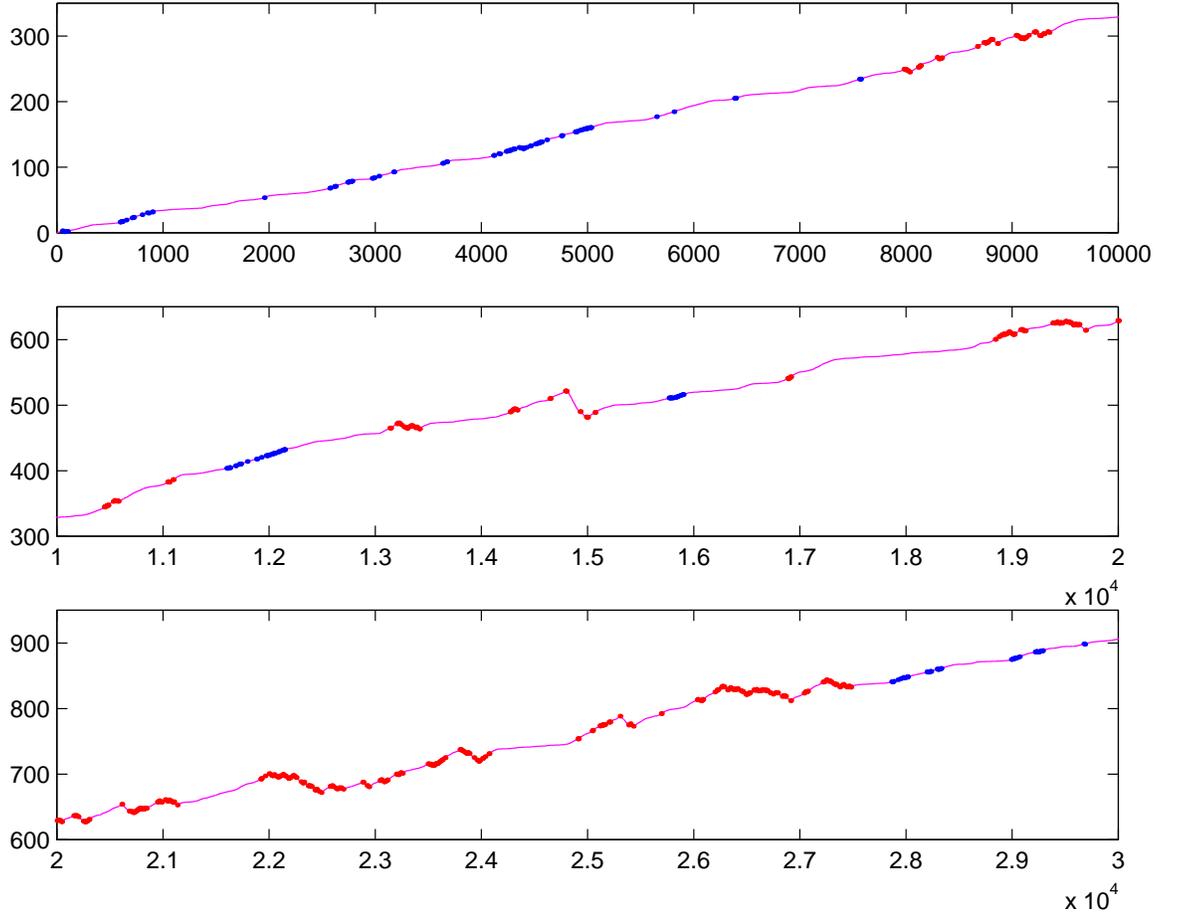}
\caption{This sequence of graphs shows the aggregate accumulation of wealth during a random sample path of length $10^5$, with $d=2, N=128, \alpha=6, \lambda=1$.  The blue dots corresponds to instances when the path passes from $g_1=43$ while the red dots corresponds to crossings of $g_4=85$.}
\label{fig:accumwealth}
\end{figure}

\begin{figure}
\epsfxsize=6in
\epsfbox{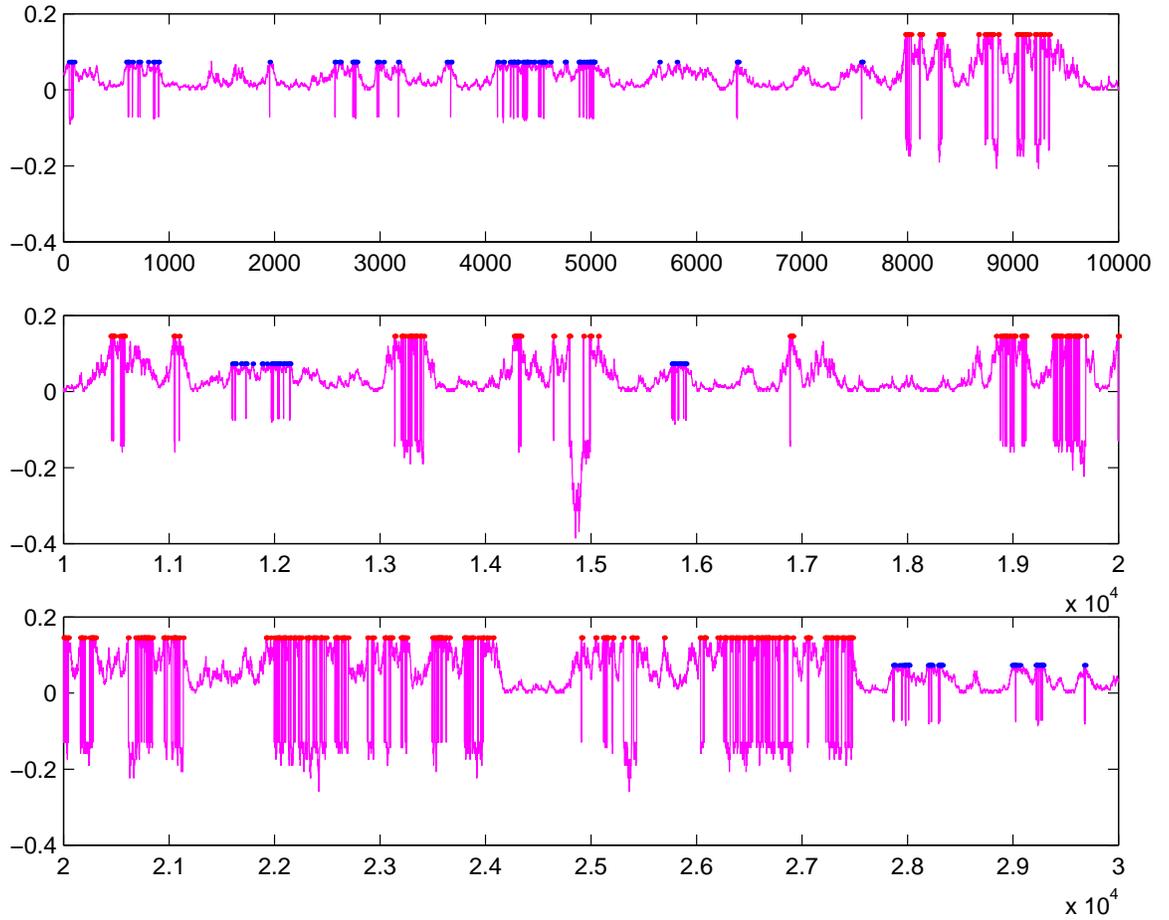}
\caption{This sequence shows the evolution of the expected increments of the aggregate wealth along the same path as Figure \ref{fig:accumwealth}.  Here the red dots populate the highest peaks, while the blue dots are restricted to the second-highest set of peaks.}
\label{fig:wealthinc}
\end{figure}

Figures \ref{fig:accumwealth} and \ref{fig:wealthinc} illustrate the importance of the crossings of $g_1$ and $g_4$ for a random path of the wealth process.  In this instance, $\lambda$ (and therefore price volatility) is low, and so the two intervals have merged, leaving only the interval $[43,85]$.  Observe that the crossings of the upper bound lead to substantially more volatile wealth than the crossings of the lower bound, which correspond to significantly smaller fluctuations of the expected wealth increments.

\section{Risk analysis}

We can extend the analysis that led to the martingale characterization of the wealth process to obtain an expression for the variance of the aggregate and individual wealth:
\begin{eqnarray}
V \left(W(T_k) -W(T_{k-1}) \left| {\mathcal B}_{k-1} \right. \right) & = & P(T_{k-1})^2 \left(e^{\frac {2 \lambda}{N}} -1 \right)^2 \nonumber \\
& & \left\{[2N^+ (T_{k-1}) -N+1]^2 P_{\scriptscriptstyle -+} \left(1 - P_{\scriptscriptstyle -+} \right) + \right. \nonumber \\
& & + e^{-{\frac {4 \lambda}{N}}} [2N^+ (T_{k-1}) -N-1 ]^2 P_{\scriptscriptstyle +-} \left(1 - P_{\scriptscriptstyle +-} \right) + \nonumber \\
& & \left. + 2 e^{-{\frac {2 \lambda}{N}}} P_{\scriptscriptstyle -+} P_{\scriptscriptstyle +-} \left([ 2N^+ (T_{k-1}) -N]^2 -1 \right) \right\} \label{eq:variance}
\end{eqnarray}

\begin{figure}
\epsfxsize=4in
\epsfbox{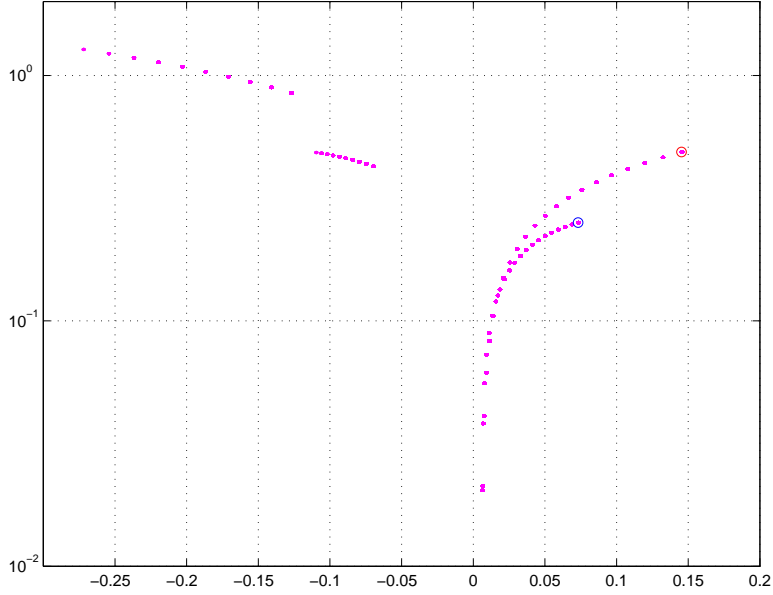}
\caption{This graph shows the standard deviation of $\Delta W$ as a function of ${\bf E} [\Delta W]$.}
\label{fig:frontier}
\end{figure}

We proceed to apply (\ref{eq:variance}) to the sample paths of the our stochastic process and investigate the relationship between expected gains/losses and the standard deviation of the wealth increments (i.e. the square root of the expression in (\ref{eq:variance}) as a measure of risk.  Figure \ref{fig:frontier} illustrates this relationship.  The graph exhibits four disconnected branches.  Two branches correspond to expected losses, for which the standard deviation has a fixed negative slope.  The branch with the higher risk corresponds to excursions of $N^+$ above $g_4$ while the lower risk branch corresponds to excursions below $g_1$.  

The two remaining branches in Figure \ref{fig:frontier} correspond to gains.  Both exhibit a significantly nonlinear, increasing relationship between risk and gain/loss.  The blue circle at the end of the lower risk branch corresponds to crossings of $g_1$.  Specifically, the entire lower right branch corresponds to fluctuations of $N^+$ between $g_1$ and $N/2$, which corresponds to the vertical asymptote towards zero risk.  When the path crosses $N/2$ and approaches $g_4$, then we find ourselves on the upper right branch, which culminates with the red circle corresponding to crossings of $g_4$.  

From this analysis we can conclude that stochastic fluctuations of the market away from the stable equilibrium with surplus buyers introduces more risk to the aggregate wealth process than those around the stable equilibrium with surplus sellers.

\section{Sojourn times of the wealth process}

\begin{figure}
\epsfxsize=5in
\epsfbox{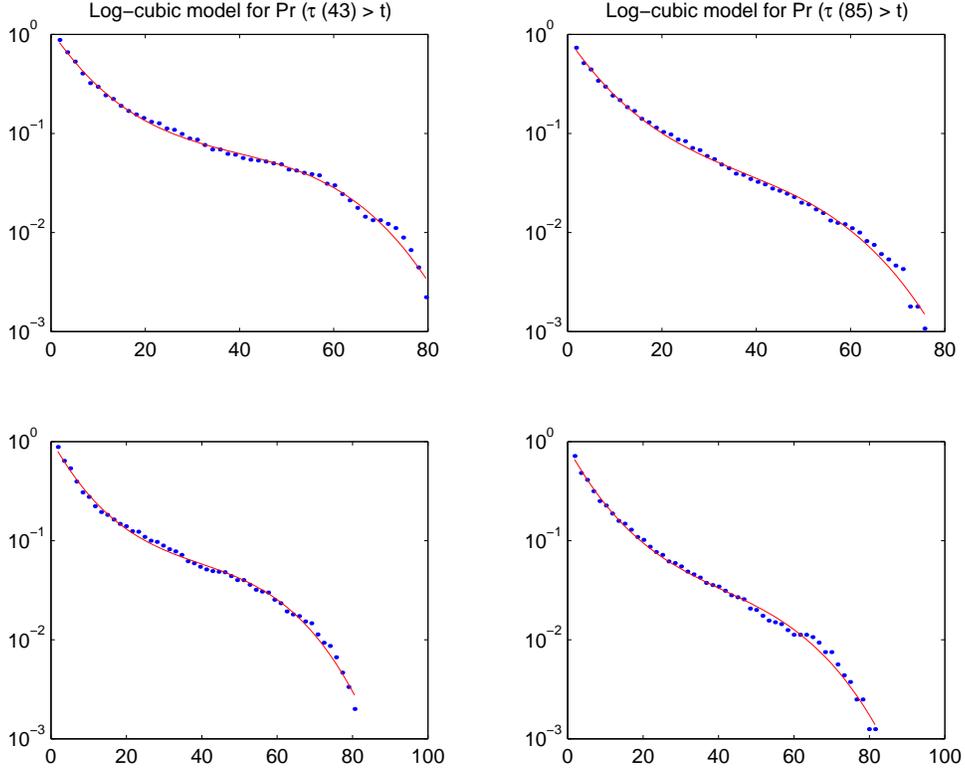}
\caption{Spline models for the tails of the inter-crossing distributions for a random path with $d=2, N=128, \alpha=6$; $\lambda$ is 1 for the top two and 5 for the bottom two graphs.}
\label{fig:tails}
\end{figure}

We complete our analysis by investigating the distribution of times the weath process spends following a sub- or super-martingale.  Let $\tau_\ell (i)$ denote the duration of the $\ell^{\rm th}$ and $(\ell+1)^{\rm st}$ crossings of the state $N^+ =i$.  Specifically, for each $i$, define recursively the infinite sequence of stopping times
$$\tilde{T}_0 (i)=0,$$
$$\tilde{T}_\ell (i)=\inf_k \left\{T_k> \tilde{T}_{\ell-1} (i) \left| N^+ (T_k) = i \mbox {  and  } \bar{X}_k \bar{X}_{k+1} >0 \right. \right\}.$$
Then
$$\tau_\ell (i) = \tilde{T}_\ell (i) - \tilde{T}_{\ell-1} (i),$$
is the stochastic process we are interested in.  Figure \ref{fig:tails} shows the tails of the $\tau_\ell$ as $\ell$ increases.  We observe that a cubic polynomial fits the logarithm of these tails admirably, leading to a model of the type
$$\lim_{\ell \rightarrow \infty} {\rm Pr} \left(\tau_\ell (i) >t \right) \approx \exp \left\{C_3 t^3 + C_2 t^2 + C_1 t +C_0 \right\}.$$
We have numerical evidence that these coefficients are robust to uncertainty about $\lambda$ and $\alpha$, so long as the qualitative behavior as described by Theorem \ref{theorem1} does not change.

\section{Conclusions and next steps}

In this paper we derived some properties of the stochastic process that describes the wealth accumulated by agents in the market microstructure model introduced in \cite{theo1}.  Specifically, our intention was to investigate a putative incentive structure which motivates the agents to act according to the interacting particle dynamics \cite{durrett}.  If successful, such an endeavor will likely help bridge the gap between mechanistic models for the emergence of randomness in market microstructure and their behavioral counterparts.

Our results so far have allowed us to take some steps in this direction.  First, we were able to deduce the stable equilibria of the stochastic system by  studying the strategic incentives that each agent faces in their attempt to secure a submartingale wealth process.  In this context, it is interesting to observe that we can use our analysis to predict that, in a market of the type modeled here, there will generally occur instances of discontinuous jumps of the equilibrium configuration, as necessitated by qualitative changes to the submartingale ranges from Theorem \ref{theorem1}.  For example, when $g_2 \rightarrow g_1$, the stable equilibrium at $g_1$ eventually disappears.  The small fluctuations around $g_1$ that were supported by its status as a stable equilibrium are no longer tenable, leading to large fluctuations that wil bring the market quickly in the neighborhood of $g_4$.  

Second, we saw that the process indeed tends to spend the majority of time in the neighborhood of the strategically preferred equilibria, leading the wealth process to approximate a martingale.  This result can be explained equally well mechanistically using the invariant measure for the Markov process \cite{theo1}, or behaviorally on the basis of the incentive structure faced by the individual agents.

Third, we obtained some insight into some inherent asymmetries in this market microstructure model.  The two stable equilibria are only symmetric for very low volatility.  Even then, the risk is different in the neighborhood of each one.  As price volatility increases, so does the asymmetry of the wealth process, passing first from a complex intermediate stage characterized by two disconnected intervals of submartingale behavior, to arrive finally at a unique equilibrium with significant buyer surplus.  

Our results have motivated us to look more closely at the path properties of the wealth process.  A more thorough analysis of the excursions away from the two equilibria may shed some light into the nature of the risk asymmetry observed above.  Furthermore, we believe that a deeper understanding of the sojourn time distributions will help us characterize the properties of the accumulated growth process, which is the stochastic integral of the increments we investigated here with respect to the invariance measure on $N^+$.

Finally, our interest lies in adding explicit strategic degrees of freedom to the agents.  One possibility that we have begun investigating is to allow agents to maintain a memory of the inputs they receive from their local neighborhood as well as the global imbalance.  In such an environment we would like to know if any agent has an incentive to increase their memory.  The analysis of the wealth process presented here can serve as the basis for such behavioral extensions of the original spin market model, by providing the incentive structure to motivate the agents to behave strategically.

\bibliographystyle{amsalpha}

\end{document}